\documentclass[12pt,letterpaper]{amsart}
\usepackage{amssymb}

\newcommand{\la}{\langle}
\newcommand{\ra}{\rangle}
\newcommand{\R}{\text{Re}\,}
\newcommand{\I}{\text{Im}\,}

\usepackage{amssymb}
\usepackage{amsthm}
\usepackage{amsxtra}
\usepackage{multirow}
\usepackage{graphicx}
\newtheorem{theorem}{Theorem}
\newtheorem{definition}[theorem]{Definition}
\newtheorem{proposition}[theorem]{Proposition}
\newtheorem{lemma}[theorem]{Lemma}

\numberwithin{equation}{section}
\numberwithin{theorem}{section}

\title[IBVP for 1D NLS on the half-line]{The initial-boundary value problem for the 1D nonlinear Schr\"odinger equation on the half-line}
\author{Justin Holmer}
\address{University of California, Berkeley}

\subjclass{35Q55}
\keywords{nonlinear Schr\"odinger equation (NLS), initial-boundary value problem, Cauchy problem, well-posedness}

\thanks{The content of this article appears as part of the author's Ph.D. thesis at the University of Chicago.}

\begin{document}

\begin{abstract}
We prove,  by adapting the method of Colliander-Kenig \cite{CK02}, local well-posedness of the initial-boundary value problem for the one-dimensional nonlinear Schr\"odinger equation $i\partial_tu +\partial_x^2u +\lambda u|u|^{\alpha-1}=0$ on the half-line under low boundary regularity assumptions.
\end{abstract}

\maketitle

\tableofcontents

\newpage
\section{Introduction}
We consider the initial-boundary value problem on the right half-line for the one-dimensional nonlinear Schr\"odinger (1D NLS) equation 
\begin{equation} \label{SE:1}
\left\{
\begin{aligned}
&i\partial_t u + \partial_x^2 u + \lambda u |u|^{\alpha-1} = 0  && \text{for }(x,t)\in (0,+\infty)\times (0,T)\\
&u(0,t) = f(t) && \text{for }t\in (0,T)\\
&u(x,0) = \phi(x) && \text{for }x\in (0,+\infty)
\end{aligned}
\right.
\end{equation}
where $\lambda\in \mathbb{C}$.  

On $\mathbb{R}$, we define the homogeneous $L^2$-based Sobolev spaces $\dot{H}^s=\dot{H}(\mathbb{R})$ by the norm $\|\phi \|_{\dot{H}^s}= \| |\xi|^s \hat{\phi}(\xi) \|_{L^2_\xi}$ and the $L^2$-based inhomogeneous Sobolev spaces $H^s=H^s(\mathbb{R})$ by the norm $\| \phi \|_{H^s} = \| \la \xi \ra^s \hat{\phi}(\xi) \|_{L^2_\xi}$, where $\la \xi \ra = (1+|\xi|^2)^{1/2}$.  In addition, we shall need $L^2$-based inhomogeneous Sobolev spaces on the half-line $\mathbb{R}^+=(0,+\infty)$, which we denote $H^s(\mathbb{R}^+)$.  These are defined, for $s\geq 0$, as: $\phi \in H^s(\mathbb{R}^+)$ if $\exists \; \tilde{\phi} \in H^s(\mathbb{R})$ such that $\tilde{\phi}(x)=\phi(x)$ for a.e.\ $x>0$; in this case we set $\|\phi \|_{H^s(\mathbb{R}^+)} = \inf_{\tilde{\phi}} \|\tilde{\phi} \|_{H^s(\mathbb{R})}$.  We also similarly define, for $s\geq 0$,  $\phi\in H^s(0,L)$ if $\exists \; \tilde{\phi}\in H^s(\mathbb{R})$ such that $\phi(x)=\tilde{\phi}(x)$ a.e.\ on $(0,L)$; in this case we set $\|\phi\|_{H^s(0,L)} = \inf_{\tilde{\phi}} \|\tilde{\phi}\|_{H^s}$. 

The \textit{local smoothing} inequality of \cite{KPV91} for the 1D Schr\"odinger group is
\begin{equation*} 
\|e^{it\partial_x^2}\phi \|_{L_x^\infty \dot{H}_t^\frac{2s+1}{4}} \leq c \|\phi\|_{\dot{H}^s}
\end{equation*}
This inequality is sharp in the sense that $\frac{2s+1}{4}$ cannot be replaced by any higher number.  We are thus motivated to consider initial-boundary data pairs $(\phi(x), f(t)) \in H^s(\mathbb{R}_x^+)\times H^\frac{2s+1}{4}(\mathbb{R}_t^+)$ and inclined to consider this configuration optimal in the scale of $L^2$-based Sobolev spaces. 

 Note that the trace map $\phi \to \phi(0)$ is well-defined on $H^s(\mathbb{R}^+)$ when $s>\frac{1}{2}$.  Thus, if $s>\frac{1}{2}$, then $\frac{2s+1}{4}>\frac{1}{2}$ and both $\phi(0)$  and $f(0)$ are well-defined quantities.  Since $\phi(0)$ and $f(0)$ are both meant to represent $u(0,0)$, they must agree.  

Therefore, we consider \eqref{SE:1} for $0\leq s<\frac{3}{2}$ in the setting
\begin{equation} \label{E:111}
\phi\in H^s(\mathbb{R}^+), \; f\in H^\frac{2s+1}{4}(\mathbb{R}^+), \; \text{and if }\tfrac{1}{2}<s<\tfrac{3}{2}, \; \phi(0)=f(0)
\end{equation}

The solutions we construct shall have the following characteristics.
\begin{definition} \label{D:strong}
$u(x,t)$ will be called a \emph{distributional solution of \eqref{SE:1}, \eqref{E:111} on $[0,T^*)$ with strong traces} if
\begin{enumerate}
\item \label{I30a}$u$ belongs to a space $X$ with the property that $u\in X$ implies $u|u|^{\alpha-1}$ is defined as a distribution.
\item \label{I30} $u(x,t)$ satisfies the equation \eqref{SE:1} in the sense of distributions on the set $(x,t) \in (0,+\infty) \times (0,T^*)$.
\item \label{I31}\emph{Space traces:} $\forall \; T<T^*$, we have $u\in C( [0,T]; \; H^s_x)$ and  $u(\cdot,0) = \phi$ in  $H^s(\mathbb{R}^+)$.
\item \emph{Time traces:} $\forall \; T<T^*$, we have $u\in C( \mathbb{R}_x ; H^\frac{2s+1}{4}(0,T))$ and $u(0, \cdot )=f$ in $H^\frac{2s+1}{4}(0,T)$.
\end{enumerate}
\end{definition}

For the purposes of uniqueness in the high regularity setting $s>\frac{1}{2}$, we can consider a weaker notion of solution.

\begin{definition}
$u(x,t)$ will be called a \emph{distributional solution of \eqref{SE:1}, \eqref{E:111} on $[0,T^*)$ with weak traces} if it satisfies conditions \ref{I30a}, \ref{I30} of Definition \ref{D:strong} and
\begin{enumerate} \addtocounter{enumi}{2}
\item \emph{One-sided space traces:} $\forall \; T<T^*$, we have $u\in C( [0,T]; \; H^s(\mathbb{R}_x^+))$ and $u(\cdot,0) = \phi$ in $H^s(\mathbb{R}^+)$.
\item \emph{Boundary values:} $\forall \; T<T^*$, we have $\displaystyle \lim_{x \downarrow 0} \| u(x,\cdot)-f \|_{H^\frac{2s+1}{4}(0,T)} = 0$.
\end{enumerate}
\end{definition}

So that we may, at a later time, properly address the matter of uniqueness in the low regularity  $s<\frac{1}{2}$ setting, we shall introduce the concept of mild solution used by \cite{BSZ04}.

\begin{definition}
$u(x,t)$ is a \emph{mild solution} of \eqref{SE:1} on $[0,T^*)$ if $\forall \; T<T^*$,  $\exists$ a sequence $\{  u_n  \}$ in $C( [0,T]; \; H^2(\mathbb{R}_x^+) )\cap C^1([0,T]; \; L^2(\mathbb{R}_x^+) )$ such that 
\begin{enumerate}
\item $u_n(x,t)$ solves \eqref{SE:1} in $L^2(\mathbb{R}_x^+)$ for $0<t<T$.
\item $\displaystyle \lim_{n\to +\infty} \|u_n -u \|_{C([0,T];\, H^s(\mathbb{R}_x^+))} =0$.
\item $\displaystyle \lim_{n\to +\infty} \|u_n(0,\cdot)-f\|_{H^\frac{2s+1}{4}(0,T)} =0$.
\end{enumerate}
\end{definition}

\cite{BSZ04} have announceed a method for proving uniqueness of mild solutions for the Korteweg-de Vries (KdV) equation on the half-line (to be discussed further in \cite{BSZ05}), and the techniques of this forthcoming paper may also apply here to resolve the uniqueness problem for $0\leq s <\frac{1}{2}$.

We establish in \S \ref{S:uniqueness} the following straightforward fact.
\begin{proposition} \label{P:uniqueness}
For $s>\frac{1}{2}$, $u$ is a distributional solution of \eqref{SE:1}, \eqref{E:111} with weak traces if and only if it is a mild solution; in this case $u$ is unique.
\end{proposition}

Our main result is the following existence statement.
\begin{theorem} \label{T:main} \quad
\begin{enumerate}
\item \label{I:mainsub} \emph{Subcritical:} Suppose 
$$\textstyle 0\leq s<\frac{1}{2}, \text{ and } 2\leq  \alpha < \frac{5-2s}{1-2s}$$ 
or 
$$\textstyle \frac{1}{2}< s<\frac{3}{2}, \text{ and } 2\leq  \alpha < \infty$$
Then $\exists \; T^*>0$ and $u$ that is both a mild solution and a distributional solution with strong traces of \eqref{SE:1},\eqref{E:111} on $[0,T^*)$.  If $T^*<\infty$, then $\lim_{t\uparrow {T^*}} \|u(\cdot,t)\|_{H_x^s} = \infty$.  Also, $\forall \; T<T^*$, $\exists \; \delta_0=\delta_0(s,T,\phi,f)>0$ such that if $0<\delta\leq \delta_0$ and $\|\phi-\phi_1\|_{H^s(\mathbb{R}^+)} + \| f-f_1\|_{H^\frac{2s+1}{4}(\mathbb{R}^+)}<\delta$ then there is a solution $u_1$(as above) on $[0,T]$, corresponding to $(\phi_1,f_1)$, such that  $\|u-u_1\|_{C(  [0,T]; \; H_x^s )} + \|u-u_1\|_{C( \mathbb{R}_x; \; H^\frac{2s+1}{4}(0,T))} \leq c\delta$, with $c=c(s,T,f,\phi)$.

\item \label{I:maincrit} \emph{Critical:} Suppose $0\leq s<\frac{1}{2}$ and $\alpha= \frac{5-2s}{1-2s}$.  Then $\exists \;T^*>0$ maximal and $u$ that is both a mild solution and a distributional solution with strong traces of \eqref{SE:1},\eqref{E:111} on $[0,T^*)$.  Also, $\exists \; T=T(s,\phi,f)<T^*$ and $\exists \; \delta_0=\delta_0(s,\phi,f)>0$ such that if $0<\delta\leq \delta_0$ and $\|\phi-\phi_1\|_{H^s(\mathbb{R}^+)} + \| f-f_1\|_{H^\frac{2s+1}{4}(\mathbb{R}^+)}<\delta$ then there is a solution $u_1$(as above) on $[0,T]$, corresponding to $(\phi_1,f_1)$, such that  $\|u-u_1\|_{C(  [0,T]; \; H_x^s )} + \|u-u_1\|_{C( \mathbb{R}_x; \; H^\frac{2s+1}{4}(0,T))} \leq c\delta$, with $c=c(s,f,\phi)$.
\end{enumerate}
\end{theorem}
\noindent Note that in \ref{I:maincrit}, we may not have blow-up in the norm $\|u(\cdot,t)\|$ as $t\uparrow T^*$.

The proof of Theorem \ref{T:main} involves the introduction of a boundary forcing operator analogous to that introduced by \cite{CK02} in their treatment of the generalized Korteweg de-Vries equation (gKdV) on the half-line, and incorporates the techniques of the standard proof of local well-posedness of the corresponding initial-value problem based on the Strichartz estimates (see \cite{CW90}).  

One could also consider the left half-line problem 
\begin{equation*} 
\left\{
\begin{aligned}
&i\partial_t u + \partial_x^2 u + \lambda u |u|^{\alpha-1} = 0  && \text{for }(x,t)\in (-\infty,0)\times (0,T)\\
&u(0,t) = f(t) && \text{for }t\in (0,T)\\
&u(x,0) = \phi(x) && \text{for }x\in (-\infty,0)
\end{aligned}
\right.
\end{equation*}
although this is actually identical to the right half-line problem \eqref{SE:1} by the transformation $u(x,t)\to u(-x,t)$.

We plan, in a future publication, to examine the initial-boundary value problem for the line-segment
\begin{equation*} 
\left\{
\begin{aligned}
&i\partial_t u + \partial_x^2 u + \lambda u |u|^{\alpha-1} = 0  && \text{for }(x,t)\in (0,L)\times (0,T)\\
&u(0,t) = f_1(t) && \text{for }t\in (0,T)\\
&u(L,t) = f_2(t) && \text{for }t\in (0,T) \\
&u(x,0) = \phi(x) && \text{for }x\in (0,L)
\end{aligned}
\right.
\end{equation*}
and consider global existence questions for the half-line and line-segment problems.

We now briefly mention some earlier work and alternate perspectives on this problem and related problems.  The main new feature of our work is the low regularity requirements for $\phi$ and $f$.  Under higher regularity assumptions, more general results are already available.  \cite{MR2002d:35196} considered a bounded or unbounded domain $\Omega\subset\mathbb{R}^n$ with smooth boundary $\partial\Omega$,  and proved global existence of solutions to
\begin{equation} \label{E:110}
\left\{\begin{aligned}
&i\partial_t u + \Delta u + \lambda u|u|^{\alpha-1} = 0 && \text{for }(x,t) \in \Omega\times (0,T)\\
&u(x,t) = f(x,t) && \text{for }x\in \partial\Omega \\
&u(x,0) = \phi(x) && \text{for }x\in \Omega
\end{aligned}
\right.
\end{equation}
where $f\in C^3(\partial\Omega)$ is compactly supported, $\phi\in H^1(\Omega)$, and $\lambda <0$.  This solution is obtained as a limit of solutions to approximate problems after several \textit{a priori} identities have been established.   Earlier, \cite{MR92d:35268} and \cite{MR2001i:35259} had obtained solutions to \eqref{SE:1} for $\alpha>3$, $\lambda<0$ and $\alpha=3$, $\lambda\in \mathbb{R}$ for $\phi \in H^2(\mathbb{R}^+)$ and $f\in C^2(0,T)$, using semigroup techniques and \textit{a priori} estimates.  The  problem \eqref{E:110} with $f=0$ had been considered previously  (\cite{MR81i:35139}  \cite{MR85f:35064}  \cite{MR92j:35175} \cite{MR90h:35224}  \cite{MR2001m:35291}).

\cite{MR2004d:37100} in the integrable case $\alpha=3$, $\lambda=\pm 2$ with  $\phi$ Schwartz and $f$ sufficiently smooth, obtained a solution to \eqref{SE:1} by reformulating the problem as a $2\times 2$ matrix Riemann-Hilbert problem.  In this setting, \cite{MR2033706} obtain an explicit representation for $\partial_x u(0,t)$. 

Outline:  In \S \ref{S:N}, we discuss some notation, introduce function spaces  and recall some needed properties of these function spaces.  In \S \ref{S:RL}, we review the definition and basic properties of the Riemann-Liouville fractional integral.  In \S \ref{S:G}, \ref{S:I}, we state the needed estimates for the group and inhomogeneous solution operator.  In \S \ref{S:BF}, we define the boundary forcing operator, adapted from \cite{CK02}, and prove the needed estimates for it.  In \S \ref{S:proofmain}, we prove Theorem \ref{T:main}.  In \S \ref{S:uniqueness}, we prove Prop.\ \ref{P:uniqueness}.

\section{Notations and some function space properties} \label{S:N}

Let $\chi_S$ denote the characteristic function for the set $S$.  We shall write $L_T^q$ to mean $L^q([0,T])$.
Set $\hat{\phi}(\xi) = \int_x e^{-ix\xi}\phi(x)\, dx$.  Define $(\tau-i0)^{\alpha}$ as the limit, in the sense of distributions, of $(\tau+i\gamma)^{-\alpha}$ as $\gamma\uparrow 0$.  Let $\la \xi \ra^{s} = (1+|\xi|^2)^{s/2}$.  Let $\widehat{D^sf}(\xi)=|\xi|^s\hat{f}(\xi)$.  The homogeneous $L^2$-based Sobolev spaces are $\dot{H}^s(\mathbb{R}) = (-\partial^2)^{-s/2}L^2(\mathbb{R})$ and the inhomogeneous $L^2$-based $H^s(\mathbb{R})=(1-\partial^2)^{-s/2}L^2(\mathbb{R})$.  We also set, for $1\leq p \leq \infty$, $W^{s,p} = (I-\partial^2)^{-s/2}L^p$.  We use the notation $H^s$ to mean $H^s(\mathbb{R})$ (and not $H^s(\mathbb{R}^+)$ or $H_0^s(\mathbb{R}^+)$). The trace operator $\phi \mapsto\phi(0)$ is defined for $\phi\in H^s(\mathbb{R})$ when $s>\frac{1}{2}$.  For $s\geq 0$, define $\phi \in H^s(\mathbb{R}^+)$ if $\exists \; \tilde{\phi} \in H^s(\mathbb{R})$ such that $\tilde{\phi}(x)=\phi(x)$ for $x>0$; in this case we set $\|\phi\|_{H^s(\mathbb{R}^+)} = \inf_{\tilde{\phi}} \|\tilde{\phi} \|_{H^s(\mathbb{R})}$.  For $s \geq 0$, define $\phi \in H_0^s(\mathbb{R}^+)$ if, when $\phi(x)$ is extended to $\tilde{\phi}(x)$ on $\mathbb{R}$ by setting $\tilde{\phi}(x)=0$ for $x<0$, then $\tilde{\phi}\in H^s(\mathbb{R})$; in this case we set $\| \phi \|_{H_0^s(\mathbb{R}^+)} = \|\tilde{\phi}\|_{H^s(\mathbb{R})}$. Define $\phi \in C_0^\infty(\mathbb{R}^+)$ if $\phi \in C^\infty(\mathbb{R})$ with $\text{supp }\phi \subset [0,+\infty)$ (so that, in particular, $\phi$ and all of its derivatives vanish at $0$), and $C_{0,c}^{\infty}(\mathbb{R}^+)$ as those members of $C_0^\infty(\mathbb{R}^+)$ with compact support.  We remark that $C_{0,c}^{\infty}(\mathbb{R}^+)$ is dense in $H_0^s(\mathbb{R}^+)$ for all $s\in \mathbb{R}$.  We shall take a fixed $\theta\in C_c^\infty(\mathbb{R})$ such that $\theta(t)=1$ on $[-1,1]$ and $\text{supp}\, \theta\subset [-2,2]$.   Denote by $\theta_T(t) = \theta(tT^{-1})$.

\begin{lemma}[\cite{CK02} Lemma 2.8]\label{CK28}  If $0\leq \alpha < \frac{1}{2}$, then $\| \theta_T h \|_{H^\alpha} \leq c\la T \ra^{\alpha}\|h\|_{\dot{H}^\alpha}$, where $c=c(\alpha, \theta)$.
\end{lemma}

\begin{lemma}[\cite{JK95} Lemma 3.5] \label{JK35}
If $-\frac{1}{2}< \alpha< \frac{1}{2}$, then $\| \chi_{(0,+\infty)}f \|_{H^\alpha} \leq c \| f \|_{H^\alpha}$, where $c=c(\alpha)$.
\end{lemma}

\begin{lemma}[\cite{CK02} Prop.\ 2.4, \cite{JK95} Lemma 3.7, 3.8] \label{JK37}
If $\frac{1}{2}<\alpha<\frac{3}{2}$, then $ H_0^\alpha(\mathbb{R}^+)=\{ f\in H^\alpha(\mathbb{R}^+) \; \mid \; f(0)=0 \}$ and if $f\in H^\alpha(\mathbb{R}^+)$ with $f(0)=0$, then $\|\chi_{(0,+\infty)} f\|_{H_0^\alpha(\mathbb{R}^+)} \leq c\|f\|_{H^\alpha(\mathbb{R}^+)}$, where $c=c(\alpha)$.
\end{lemma}

The following Gronwall-type inequality can be obtained by applying the H\"older inequality iteratively:

\begin{lemma} \label{L:Gronwall}
If $1\leq q_1<q\leq \infty$ and $\forall \; t\geq 0$ 
$$\left( \int_0^t |g(s)|^q \, ds \right)^{1/q} \leq c \delta + c\left( \int_0^t |f(s)|^{q_1} \, ds \right)^{1/{q_1}}$$
then with $\gamma$ defined by $2c\gamma^{\frac{1}{q_1}-\frac{1}{q}} = 1$, we have $\forall \; t\geq 0$,
$$\left(\int_0^t |f(s)|^{q_1} \, ds \right)^{1/{q_1}} \leq (\gamma t)^{\gamma t} \delta$$

\end{lemma}

A version of the chain rule for fractional derivatives is
\begin{lemma}[Prop.\ 3.1 in \cite{CW91}] \label{CWchain}
Let $0<s<1$, $u:\mathbb{R} \to \mathbb{R}^2$ and $F:\mathbb{R}^2 \to \mathbb{R}^2$, $F\in C^1$, so that $F'(u)$ is a $2\times 2$ matrix.  Then
$$\| D^sF(u) \|_{L^r} \leq c\|F'(u)\|_{L^{r_1}} \|D^su\|_{L^{r_2}}$$
for $\frac{1}{r}=\frac{1}{r_1}+\frac{1}{r_2}$ with $1< r,r_1,r_2<\infty$.
\end{lemma}
The product rule for fractional derivatives is
\begin{lemma}[Prop.\ 3.3 in \cite{CW91}] \label{CWproduct}
Let $0<s<1$.  If $u,v: \mathbb{R}\to\mathbb{R}$, then
$$\|D^s(uv)\|_{L^r} \leq \|D^su\|_{L^{r_1}} \| v \|_{L^{r_2}} + \|u\|_{L^{r_3}}\|D^sv\|_{L^{r_4}}$$
for $1<r,r_1,r_2,r_3,r_4<\infty$ and $\frac{1}{r}=\frac{1}{r_1}+\frac{1}{r_2}$, $\frac{1}{r}=\frac{1}{r_3}+\frac{1}{r_4}$.
\end{lemma}

\section{The Riemann-Liouville fractional integral} \label{S:RL}

The tempered distribution $\frac{t_+^{\alpha-1}}{\Gamma(\alpha)}$ is defined as a locally integrable function for $\text{Re }\alpha>0$, i.e.\
$$\left< \frac{t_+^{\alpha-1}}{\Gamma(\alpha)}, f \right> = \frac{1}{\Gamma(\alpha)}\int_0^{+\infty} t^{\alpha-1} f(t) \, dt$$
Integration by parts gives, for \text{Re }$\alpha>0$, that 
$$\frac{t_+^{\alpha-1}}{\Gamma(\alpha)} = \partial_t^k \left[ \frac{t_+^{\alpha+k-1}}{\Gamma(\alpha+k)} \right]$$
for all $k\in \mathbb{N}$.
  This formula can be used to extend the definition (in the sense of distributions) of $\frac{t_+^{\alpha-1}}{\Gamma(\alpha)}$ to all $\alpha \in \mathbb{C}$.  In particular, we obtain
$$\left. \frac{t_+^{\alpha-1}}{\Gamma(\alpha)}\right|_{\alpha=0} = \delta_0(t)$$
A change of contour calculation shows that 
$$\left[\frac{t_+^{\alpha-1}}{\Gamma(\alpha)} \right]\sphat(\tau)=e^{-\frac{1}{2}\pi i \alpha}(\tau - i0)^{-\alpha}$$
where $(\tau-i0)^{-\alpha}$ is the distributional limit. If $f\in C_0^\infty(\mathbb{R}^+)$, we define
$$\mathcal{I}_\alpha f  = \frac{t_+^{\alpha-1}}{\Gamma(\alpha)} * f$$
Thus, when $\text{Re }\alpha>0$,
$$\mathcal{I}_\alpha f(t)  = \frac{1}{\Gamma(\alpha)}\int_0^t (t-s)^{\alpha-1} f(s) \, ds$$
and $\mathcal{I}_0f=f$, $\mathcal{I}_1f(t)=\int_0^t f(s) \, ds$, and $\mathcal{I}_{-1}f=f'$.  Also $\mathcal{I}_{\alpha}\mathcal{I}_\beta = \mathcal{I}_{\alpha+\beta}$, which follows from the Fourier transform formula.  For further details on the distribution $\frac{t_+^{\alpha-1}}{\Gamma(\alpha)}$, see \cite{F98}.

\begin{lemma} \label{L:RL}
If $h\in C_0^{\infty}(\mathbb{R}^+)$, then $\mathcal{I}_{\alpha}h \in C_0^{\infty}(\mathbb{R}^+)$, for all $\alpha\in \mathbb{C}$.
\end{lemma}

\begin{lemma}[\cite{H04}] \label{L:RL1}
If $0\leq \alpha < +\infty$ and $s\in \mathbb{R}$, then
\begin{equation*} 
\| \mathcal{I}_{-\alpha}h \|_{H_0^s(\mathbb{R}^+)}\leq c \|h\|_{H_0^{s+\alpha}(\mathbb{R}^+)} 
\end{equation*}
\end{lemma}

\begin{lemma}[\cite{H04}] \label{L:RL2}
If $0\leq \alpha < +\infty$, $s\in \mathbb{R}$, $\mu \in C_0^\infty(\mathbb{R})$
\begin{equation*}
\| \mu \mathcal{I}_\alpha h \|_{H_0^s(\mathbb{R}^+)} \leq c \|h\|_{H_0^{s-\alpha}(\mathbb{R}^+)}
\end{equation*}
where $c=c(\mu)$.
\end{lemma}

\section{Estimates for the group} \label{S:G}
Set
\begin{equation} \label{E:70}
e^{it\partial_x^2}\phi(x) = \frac{1}{2\pi} \int_{\xi} e^{ix\xi} e^{-it\xi^2} \hat{\phi}(\xi) \, d\xi
\end{equation}
so that 
$$\left\{\begin{aligned}
& (i\partial_t + \partial_x^2)e^{it\partial_x^2}\phi = 0  && \text{for }(x,t) \in \mathbb{R}\times \mathbb{R}\\
& e^{it\partial_x^2}\phi(x) \big|_{t=0}=\phi(x) && \text{for }x\in \mathbb{R}
\end{aligned} \right.$$
\begin{lemma} \label{L:G}
Let $s\in \mathbb{R}$.  If $\phi \in H^s(\mathbb{R})$, then
\begin{enumerate}
\item \label{I:Gs} \emph{Space traces:} $\| e^{it\partial_x^2} \phi(x) \|_{C( \mathbb{R}_t;\, H_x^s)} \leq c \|\phi \|_{H^s}$.
\item \label{I:Gt}  \emph{Time traces:}  $\| \theta_T(t) e^{it\partial_x^2}\phi(x) \|_{C(\mathbb{R}_x; \, H_t^\frac{2s+1}{4})} \leq c\la T \ra^{1/4} \|\phi\|_{H^s}$.
\item \label{I:Gm} \emph{Mixed-norm:} If $2\leq q, \, r\leq \infty$ and $\frac{1}{q}+\frac{1}{2r}=\frac{1}{4}$, then $\|e^{it\partial_x^2}\phi(x) \|_{L_t^q W_x^{s,r}} \leq c\|\phi\|_{H^s}$.
\end{enumerate}
\end{lemma}
\begin{proof}
\ref{I:Gs} is clear from \eqref{E:70}.  \ref{I:Gt} was obtained in \cite{KPV91}.  \ref{I:Gm} was obtained by \cite{S77} (see also \cite{KT98}).
\end{proof}

\section{Estimates for the Duhamel inhomogeneous solution operator} \label{S:I}

Let
$$\mathcal{D}w(x,t) = -i\int_0^t e^{i(t-t')\partial_x^2}w(x,t') \, dt'$$
Then
$$\left\{ \begin{aligned}
& (i\partial_t + \partial_x^2)\mathcal{D}w(x,t) = w(x,t) && \text{for }(x,t) \in \mathbb{R}\times \mathbb{R}\\
& \mathcal{D}w(x,0)=0 && \text{for }x\in \mathbb{R}
\end{aligned}\right.$$
\begin{lemma}  \label{L:I}
 Suppose $2\leq q, \, r\leq \infty$ and $\frac{1}{q}+\frac{1}{2r}=\frac{1}{4}$, then 
\begin{enumerate}
\item \label{I:Is} \emph{Space traces:} If $s\in \mathbb{R}$, then $\|\mathcal{D}w \|_{C( \mathbb{R}_t; H_x^s )} \leq c \|w\|_{L_t^{q'} W_x^{s,r'}}$.
\item \label{I:It} \emph{Time traces:} If $-\frac{3}{2}<s<\frac{1}{2}$, then $\| \theta_T(t)\mathcal{D}w(x,t) \|_{C(\mathbb{R}_x; \, H_t^\frac{2s+1}{4})} \leq c \la T \ra^{1/4}\| w\|_{L_t^{q'}W_x^{s,r'}}$.
\item \label{I:Im} \emph{Mixed-norm:} If $s\in \mathbb{R}$, then $\|\mathcal{D}w \|_{L_t^q W_x^{s,r}} \leq c \|w\|_{L_t^{q'}W_x^{s,r'}}$.
\end{enumerate}
\end{lemma}
\begin{proof}
\ref{I:Is} and \ref{I:Im} are due to \cite{S77} (see also \cite{KT98}).  We now prove \ref{I:It}, following the techniques of Theorem 2.3 in \cite{KPV93}.  We use the representation
\begin{align*}
\mathcal{D}w(x,t) &= 
\begin{aligned}[t]
&-\frac{i}{2}\int_{-\infty}^{+\infty} (\text{sgn}\, t')e^{i(t-t')\partial_x^2} w(x,t') \, dt'  \\
&+\frac{1}{2\pi i} \int_\tau e^{it\tau} \left[ \lim_{\epsilon \to 0^+} \frac{1}{2\pi} \int_{|\tau+\xi^2|>\epsilon} e^{ix\xi} \frac{\hat{w}(\xi,\tau)}{\tau+\xi^2} \, d\xi \right] \, d\tau
\end{aligned} \\
&=\text{I}+\text{II}
\end{align*}
and Term II can also be written
$$\text{II} \, = \frac{1}{2\pi} \int_\tau  e^{it\tau} [m(\cdot, \tau) * \hat{w}^t(\cdot,\tau)](x) \, d\tau$$
where $\hat{w}^t(\cdot,\tau)$ denotes the Fourier transform of $w(x,t)$ in the $t$-variable alone and
$$m(x,\tau) = -\frac{1}{2}\chi_{(0,+\infty)}(\tau)\frac{\exp(-|x||\tau|^{1/2})}{|\tau|^{1/2}} +\frac{1}{2}\chi_{(-\infty,0)}(\tau) \frac{\sin(|x||\tau|^{1/2})}{|\tau|^{1/2}}$$
First we treat Term I for all $s$ and all admissible pairs $q,r$.  Pairing Term I with $f(x,t)$ such that $\|f\|_{L_x^1H_t^{-\frac{2s+1}{4}}}\leq 1$, we are left to show that
$$\left\| \int_{t'} (\text{sgn}\, t') e^{-it'\partial_x^2}w(x,t')\, dt' \right\|_{H_x^s} \leq c\|w\|_{L_t^{q'}W_x^{s,r'}}$$
and
$$\left\| \int_t \theta_T(t) e^{-it\partial_x^2}f(x,t) \, dt \right\|_{H_x^{-s}}\leq c\|f\|_{L_x^1H_t^{-\frac{2s+1}{4}}}$$
The first of these follows from the proof of \ref{I:Is}, while the second is obtained by duality and Lemma \ref{L:G}\ref{I:Gt}.
We address Term II separately for $r'=2$, $q'=1$, and $r'=1$, $q'=\frac{4}{3}$; the intermediate cases follow by interpolation.  For the case $r'=2$, $q'=1$, we use the first representation of Term II with Lemma \ref{CK28}, the change of variable $\eta=-\xi^2$, and $L^2$-boundedness of the Hilbert transform on $A_2$-weighted spaces, to obtain
$$\| \theta_T(t) (\text{Term II}) \|_{H_t^\frac{2s+1}{4}} 
\begin{aligned}[t]
&\leq c\left( \int_\xi |\xi|^s |\hat{w}(\xi,\tau)|^2 \, d\xi \right)^{1/2} \\
&\leq c\left( \int_\xi |\xi|^s \left( \int_t|\hat{w}^x(\xi,t)| \,dt \right)^2 d\xi \right)^{1/2}
\end{aligned}$$
where $\hat{w}^x(\xi,t)$ denotes the Fourier transform in the $x$-variable alone.  Complete the bound by applying Minkowskii's integral inequality and the Placherel theorem.  The validity of this step is restricted to $-\frac{3}{2}<s<\frac{1}{2}$.  

We shall only prove the $r'=1$, $q'=\frac{4}{3}$ case for $s=0$.  Note that by the second representation for Term II, $\| (\text{Term II}) \|_{L_x^\infty H_t^{1/4}}$ is 
$$ \int_\tau \int_y |\tau|^{-1/2}m(x-y,\tau)\hat{w}^t(y,\tau) \, dy \; \overline{\int_z |\tau|^{-1/2}m(x-z,\tau)\hat{w}^t(z,\tau) \, dz} \; d\tau$$
which is equivalent to
$$\int_{y,s,z,t} K(y,s,z,t) w(y,s) \overline{w(z,t)} \, dy\, ds \, dz \, dt$$
where
$$K(y,s,z,t) = \int_\tau |\tau|^{1/2}e^{-i(s-t)\tau} m(x-y,\tau) \overline{m(x-z,\tau)} \, d\tau$$
From the definition of $m$, we see that $|K(y,s,z,t)| \leq c|s-t|^{-1/2}$.  We conclude by applying the theorem on fractional integration (see Theorem 1 of Chapter V in \cite{St70}).
\end{proof}

\section{Estimates for the Duhamel boundary forcing operator}\label{S:BF}

For $f\in C_0^\infty(\mathbb{R}^+)$, define the boundary forcing operator
\begin{align}
\mathcal{L}f(x,t)&= 2e^{i\frac{1}{4}\pi} \int_0^t e^{i(t-t')\partial_x^2}\delta_0(x)\mathcal{I}_{-1/2}f(t') \, dt' \label{E:100}\\
&=\frac{1}{\sqrt{\pi}} \int_0^t (t-t')^{-1/2} \exp \left( \frac{ix^2}{4(t-t')} \right) \mathcal{I}_{-1/2}f(t') \, dt' \label{E:101}
\end{align}
The equivalence of the two definitions is evident from the formula
$$\left[ \frac{e^{-i\frac{\pi}{4}\text{sgn}\,t}}{2\sqrt{\pi}}  \frac{1}{|t|^{1/2}} \exp \left( \frac{ix^2}{4t} \right) \right]\sphat(\xi) = e^{-it\xi^2}$$
From these two definitions, we see that
$$\left\{
\begin{aligned}
&(i\partial_t + \partial_x^2)\mathcal{L}f(x,t) = 2e^{i\frac{3}{4}\pi} \delta_0(x) \mathcal{I}_{-1/2}f(t) && \text{for }(x,t)\in \mathbb{R}\times \mathbb{R} \\
&\mathcal{L}f(x,0) = 0  && \text{for }x\in \mathbb{R}\\
&\mathcal{L}f(0,t)=f(t) && \text{for }t\in \mathbb{R}
\end{aligned}
\right.$$
We now establish some continuity properties of $\mathcal{L}f(x,t)$ when $f$ is suitably nice.
\begin{lemma}
Let $f\in C_{0,c}^\infty(\mathbb{R}^+)$.
\begin{enumerate}
\item For fixed $t$, $\mathcal{L}f(x,t)$ is continuous in $x$ for all $x\in \mathbb{R}$ and  $\partial_x \mathcal{L}f(x,t)$ is continuous in $x$ for $x\neq 0$ with 
\begin{equation} \label{E:150}
\lim_{x\uparrow 0} \partial_x\mathcal{L}f(x,t) = e^{-\frac{1}{4}\pi i}\mathcal{I}_{-1/2}f(t) \qquad \lim_{x\downarrow 0} \partial_x\mathcal{L}f(x,t) = -e^{-\frac{1}{4}\pi i}\mathcal{I}_{-1/2}f(t)
\end{equation}
\item $\forall \; k=0,1,2,\ldots$ and for fixed $x$, $\partial_t^k\mathcal{L}f(x,t)$ is continuous in $t$ for all $t\in \mathbb{R}$. 
\end{enumerate}
We also have the pointwise estimates, for $k=0,1,2,\ldots$, on $[0,T]$,
$$|\partial_t^k\mathcal{L}f(x,t)| + |\partial_x \mathcal{L}f(x,t)| \leq c\la x \ra^{-N}$$
where $c=c(f,N,k,T)$.
\end{lemma}

\begin{proof}
Let us denote ``integration by parts'' by IBP.  It is clear from \eqref{E:101} and dominated convergence that, for fixed $t$, $\mathcal{L}f(x,t)$ is continuous in $x$, and for fixed $x$, $\mathcal{L}f(x,t)$ is continuous in $t$.  Let $h=2e^{i\frac{1}{4}\pi}\mathcal{I}_{-1/2}f\in C_0^{\infty}(\mathbb{R}^+)$ (by Lemma \ref{L:RL}) and $\phi(\xi,t) =  \int_0^t e^{-i(t-t')\xi} h(t') \, dt'$.  By IBP in $t'$,  $|\partial_\xi^k \phi(\xi,t)| \leq c\la \xi\ra^{-k-1}$, where $c=c(h,k,T)$, and thus 
\begin{equation}
\label{BE:926}
|\partial_\xi^k \phi(\xi^2,t)| \leq c\la \xi \ra^{-k-2}
\end{equation}
We have
\begin{equation} \label{BE:927}
\mathcal{L}f(x,t)=\int_\xi e^{ix\xi} \phi(\xi^2,t) \, d\xi
\end{equation}
and by IBP in $\xi$ and \eqref{BE:926}, we have $|\mathcal{L}f(x,t)| \leq c\la x \ra^{-N}$.  By $\partial_t[e^{i(t-t')\partial_x^2}\delta_0(x)] = -\partial_{t'}[e^{i(t-t')\partial_x^2}\delta_0(x)]$ and IBP in $t'$ in \eqref{E:100}, $\partial_t\mathcal{L}f=\mathcal{L}\partial_tf$, and thus, for fixed $x$, $\partial_t^k \mathcal{L}f(x,t)$ is continuous in $t$ and $|\partial_t^k\mathcal{L}f(x,t)|\leq c\la x \ra^{-N}$.  By $\partial_x^2[e^{i(t-t')\partial_x^2}\delta_0(x)] = i\partial_{t'}[e^{i(t-t')\partial_x^2}\delta_0(x)]$ and IBP in $t'$ in \eqref{E:100}, $\partial_x^2\mathcal{L}f(x,t)= 2e^{i\frac{3}{4}\pi} \delta_0(x)\mathcal{I}_{-1/2}f(t)-i\mathcal{L}(\partial_tf)(x,t)$.  Hence 
$$\partial_x\mathcal{L}f(x,t) = e^{i\frac{3}{4}\pi}(\text{sgn}\, x) \mathcal{I}_{-1/2}f(t) -i\int_{x'=0}^x \mathcal{L}(\partial_tf)(x',t)\, dx' + c(t)$$
Since all terms except $c(t)$ are odd in $x$, we must have $c(t)=0$.  From this we obtain \eqref{E:150}, and the bound $|\partial_x\mathcal{L}f(x,t)|\leq c$.
From \eqref{BE:927}, IBP in $\xi$ and \eqref{BE:926}, we obtain that $|\partial_x\mathcal{L}f(x,t)| \leq c|x|^{-N}$.  Combining the two previous bounds, we have $|\partial_x \mathcal{L}f(x,t)| \leq c\la x \ra^{-N}$.
\end{proof}

Now we provide an alternate representation of $\mathcal{L}f(x,t)$.
\begin{lemma} \label{L:BFalt}
 Suppose $f\in C_{0,c}^\infty(\mathbb{R}^+)$.  Then
\begin{equation}
\label{E:1}
\mathcal{L}f(x,t) = \frac{1}{2\pi} \int_\tau e^{it\tau} e^{-|x|(\tau-i0)^{1/2}} \hat{f}(\tau) \, d\tau
\end{equation}
where
$$(\tau-i0)^{\frac{1}{2}} = \chi_{(0,+\infty)}(\tau)|\tau|^{1/2} -i \chi_{(-\infty,0)}(\tau)|\tau|^{1/2}$$
\end{lemma}
\begin{proof}
It suffices to verify that 
\begin{enumerate}
\item \label{I1} On $[0,T]$, $ |\mathcal{L}f(x,t)|+ |\partial_t \mathcal{L}f(x,t)| \leq c\la x \ra^{-N}$, with $c=c(f,N,T)$.
\item \label{I2} $\mathcal{L}f(x,0) = 0$
\item \label{I3} $(i\partial_t + \partial_x^2)\mathcal{L}f(x,t) = 2\delta_0(x)e^{\frac{3}{4}\pi i}\mathcal{I}_{-1/2}f(t)$
\end{enumerate}
\ref{I1} is integration by parts in $\tau$ in \eqref{E:1} using $-2(\tau-i0)^{1/2}|x|^{-1}\partial_{\tau}[e^{-|x|(\tau-i0)^{1/2}}]=e^{-|x|(\tau-i0)^{1/2}}$.  To show \ref{I2}, note that since $f\in C_{0,c}^\infty(\mathbb{R}^+)$, $\hat{f}(\tau)$ extends to an analytic function on $\text{Im }\tau <0$ satisfying $|\hat{f}(\tau)|\leq c\la \tau \ra^{-k}$ with $c=c(f,k)$, and thus 
\begin{equation}
\label{E:2}
\mathcal{L}f(x,0) = \frac{1}{2\pi} \lim_{\gamma \uparrow 0} \int_{\I \tau = \gamma}  e^{-|x|\tau^{1/2}} \hat{f}(\tau) \, d\tau
\end{equation}
Since $|e^{-|x|\tau^{1/2}}| \leq 1$ for $\I \tau <0$, by Cauchy's theorem, $\eqref{E:2}=0$.  \ref{I3} is a direct computation from \eqref{E:1}.  

Denote the operator defined by \eqref{E:1} as $\mathcal{L}_2f(x,t)$ and the one given by \eqref{E:100}-\eqref{E:101} as $\mathcal{L}_1f(x,t)$.  Setting $w=\mathcal{L}_1f-\mathcal{L}_2f$, we have $w(x,0)=0$ and $(i\partial_t + \partial_x^2)w=0$.  Compute $\partial_t \int_x |w|^2 dx =0$, which yields $w=0$, to complete the proof.
\end{proof}

\begin{lemma} \label{L:BF}
Suppose $q,r\geq 2$ and $\frac{1}{q}+\frac{1}{2r}=\frac{1}{4}$.  
\begin{enumerate}
\item \label{I:BFs}\emph{Space traces:} If $-\frac{1}{2}<s<\frac{3}{2}$, then $\|\theta_T(t)\mathcal{L}f(x,t) \|_{C(\mathbb{R}_t; H_x^s)} \leq c \la T \ra^{1/4}\|f\|_{H_0^\frac{2s+1}{4}(\mathbb{R}^+)}$.
\item \label{I:BFt} \emph{Time traces:} If $s\in \mathbb{R}$, then $\| \mathcal{L}f \|_{C(\mathbb{R}_x; H_0^\frac{2s+1}{4}(\mathbb{R}^+_t))} \leq c \|f\|_{H_0^\frac{2s+1}{4}(\mathbb{R}^+)}$.
\item \label{I:BFm} \emph{Mixed-norm:} If $0\leq s \leq 1$, $r\neq \infty$,  we have $\|\mathcal{L}f\|_{L_t^qW_x^{s,r}} \leq c\|f\|_{H_0^\frac{2s+1}{4}(\mathbb{R}^+)}$.
\end{enumerate}
\end{lemma}
\begin{proof}
By density, it suffices to establish these facts for $f\in C_{0,c}^\infty(\mathbb{R}^+)$.  

By pairing \ref{I:BFs} with $\phi(x)$ such that $\|\phi\|_{H^{-s}}\leq 1$, we see that it suffices to show
$$\int_{t'=0}^t f(t') \theta_T(t)e^{i(t-t')\partial_x^2} \phi\big|_{x=0} \, dt' \leq c\la T \ra^{1/4}\|f\|_{H^\frac{2s+1}{4}}$$
But
$$\text{LHS} \leq \|\chi_{(-\infty,t)}f(t')\|_{H^\frac{2s+1}{4}_{t'}} \|\theta_T(t)e^{i(t-t')\partial_x^2}\phi(x) \|_{H^\frac{-2s-1}{4}_{t'}} \leq \text{RHS}$$
by Lemmas \ref{L:G}\ref{I:Gt} and \ref{JK35}.
To establish the continuity statement, write $\theta_T(t_2)\mathcal{L}f(x,t_2)-\theta_T(t_1)\mathcal{L}f(x,t_1) = \int_{t_1}^{t_2} \partial_t[\theta(t)\mathcal{L}f(x,t)]\, dt$.  By $\partial_t \mathcal{L} = \mathcal{L}\partial_t$ and the bound just derived, we have $\| \theta_T(t_2)\mathcal{L}f(x,t_2)-\theta_T(t_1)\mathcal{L}f(x,t_1) \| \leq c|t_2-t_1| \|f\|_{H_0^\frac{2s+5}{4}}$.

\ref{I:BFt} is immediate from Lemma \ref{L:BFalt}, except that we should confirm that (under the assumption $f\in C_{0,c}^\infty(\mathbb{R}^+)$, that $\partial_t^k \mathcal{L}f(x,0) =0$ for all $k=0,1,2,\ldots$.  This, however, follows from $\partial_t \mathcal{L} = \mathcal{L}\partial_t$.  The continuity statement follows by using $\mathcal{L}f(x_2,t)-\mathcal{L}f(x_1,t)= \int_{x_1}^{x_2}\partial_x \mathcal{L}f(x,t)\, dx$.  From Lemma \ref{L:BFalt}, we have
$$\partial_x \mathcal{L}f(x,t) =  e^{-\frac{1}{4}\pi i}(\text{sgn}\,x) \frac{1}{2\pi} \int_\tau e^{it\tau} e^{-|x|(\tau-i0)^{1/2}} [\mathcal{I}_{-1/2}f]\sphat\,(\tau) \, d\tau$$
and thus $$\|\mathcal{L}f(x_2,t)-\mathcal{L}f(x_1,t)\|_{H_0^\frac{2s+1}{4}(\mathbb{R}^+)} \leq c|x_2-x_1| \|f\|_{H_0^\frac{2s+3}{4}}$$

To prove \ref{I:BFm}, it suffices to establish 
\begin{equation} \label{E:162}
\| \mathcal{L}f(x,t) \|_{L_t^4 L_x^\infty} \leq c\|f\|_{\dot{H}^{1/4}}
\end{equation}
and
\begin{equation}\label{E:163}
 \| \partial_x \mathcal{L}f(x,t) \|_{L_t^4 L_x^\infty} \leq c\|f\|_{\dot{H}^{3/4}}
\end{equation}
Indeed, the proof of \ref{I:BFs} in the case $s=1$ shows
$$\|\mathcal{L}f(x,t)\|_{L_t^\infty L_x^2} \leq c \|f\|_{\dot{H}^{1/4}} \qquad \| \partial_x \mathcal{L}f(x,t) \|_{L_t^\infty L_x^2} \leq c\|f\|_{\dot{H}^{3/4}}$$
Interpolate \eqref{E:162} with the first inequality and \eqref{E:163} with the second inequality to obtain
$$\|\mathcal{L}f(x,t)\|_{L_t^q L_x^r} \leq c \|f\|_{\dot{H}^{1/4}} \qquad \| \partial_x \mathcal{L}f(x,t) \|_{L_t^q L_x^r} \leq c\|f\|_{\dot{H}^{3/4}}$$
for admissible $q,r$.  This implies
$$\| \mathcal{L}f(x,t)\|_{L_t^qW_x^{s,r}} \leq c\|f\|_{H^\frac{2s+1}{4}}, \quad r\neq \infty$$
for $s=0$ and $s=1$.  Now interpolate over $s$ between these two endpoints to obtain the result as stated.  

 By pairing LHS of \eqref{E:162} against $w(x,t) \in L_t^{4/3}L_x^1$, we see that it suffices to show
$$ \left\| \int_x\int_t e^{it\tau} e^{-|x|(\tau-i0)^{1/2}}w(x,t) \, dx\, dt \right\|_{\dot{H}^{-1/4}}$$
Writing out the $L_\tau^2$ norm, we see that it suffices to show
$$\int_{x,t,y,s} K(x,t,y,s) \, w(x,t) \overline{w(y,s)} \, dx \, dt \, dy \, ds \leq c\|w\|_{L_t^{4/3}L_x^1}$$
where
$$K(x,t,y,s) = \int_\tau |\tau|^{-1/2} e^{i(t-s)\tau} e^{-|x|(\tau-i0)^{1/2}} e^{-|y|(\tau+i0)^{1/2}} \, d\tau$$
By a change of contour calculation, it follows that $|K(x,y,t,s)| \leq c|t-s|^{-1/2}$, and hence \eqref{E:162} follows by the theorem on fractional integration.  For \eqref{E:163}, the kernel is instead
$$K(x,t,y,s) = (\text{sgn}\, x)( \text{sgn}\, y) \int_\tau |\tau|^{-1/2} e^{i(t-s)\tau} e^{-|x|(\tau-i0)^{1/2}} e^{-|y|(\tau+i0)^{1/2}} \, d\tau$$
and hence the estimation of $|K|$ is identical.
\end{proof}

\section{Existence: Proof of Theorem \ref{T:main}} \label{S:proofmain}

First we prove the subcritical assertion (a) in the case $0\leq s<\frac{1}{2}$.  Select an extension $\tilde{\phi}\in H^s$ of $\phi$ such that $\|\tilde{\phi}\|_{H^s} \leq 2\|\phi\|_{H^s(\mathbb{R}^+)}$.  Set $r=\frac{\alpha+1}{1+(\alpha-1)s}$ and $q=\frac{4(\alpha+1)}{(\alpha-1)(1-2s)}$.  This is an admissible pair with $r\geq 2$ and $q\geq 2(\frac{2}{1-2s}+1)$.  Set
$$Z = C(\mathbb{R}_t; \; H^s_x) \cap C( \mathbb{R}_x; H^\frac{2s+1}{4}_t) \cap L_t^qW_x^{s,r}$$
Take $w\in Z$.  By the chain rule (Lemma \ref{CWchain}), for $\alpha\geq 1$ (see below for details)
\begin{equation}\label{E:140}
\| D^s( |w|^{\alpha-1} w) \|_{L_{4T}^{q'} L_x^{r'}} \leq   cT^\sigma\|w\|_{L_{4T}^qW_x^{r,s}}^\alpha 
\end{equation}
for some $\sigma>0$.  Note that by Lemmas \ref{L:G}\ref{I:Gt}, \ref{L:I}\ref{I:It}, \ref{JK35}, if $w\in Z$, then $f(t)- \theta_{2T}(t)e^{it\partial_x^2}\tilde{\phi}\big|_{x=0} \in H^\frac{2s+1}{4}_0(\mathbb{R}_t^+)$ and $\theta_{2T}(t)\mathcal{D}(w|w|^{\alpha-1})(0,t) \in H^\frac{2s+1}{4}_0(\mathbb{R}_t^+)$, and the evaluation at $x=0$ in these statements is understood in the sense of $C( \mathbb{R}_x; H^\frac{2s+1}{4}_t)$.
Let
\begin{equation} \label{E:143}
\Lambda w(t) = 
\begin{aligned}[t]
&\theta_T(t)e^{it\partial_x^2}\tilde{\phi} + \theta_T(t)\mathcal{L}(f- \theta_{2T}e^{i\cdot\partial_x^2}\tilde{\phi}\big|_{x=0})(t) \\
&- \lambda \theta_T(t)\mathcal{D}(w|w|^{\alpha-1})(t) + \lambda \theta_T(t) \mathcal{L}(\theta_{2T}\mathcal{D}(w|w|^{\alpha-1})\big|_{x=0} )(t) 
\end{aligned}
\end{equation}
so that, on $[0,T]$,  $(i\partial_t + \partial_x^2)\Lambda w = - \lambda w|w|^{\alpha-1}$ for $x\neq 0$ in the sense of distributions.  By Lemmas \ref{L:G}, \ref{L:I}, \ref{L:BF} and \eqref{E:140},
\begin{equation} \label{E:141}
\|\Lambda w \|_{Z} \leq c\|\phi\|_{H^s(\mathbb{R}^+)} + c\|f\|_{H^\frac{2s+1}{4}(\mathbb{R}^+)} + cT^\sigma \|w\|_{Z}^\alpha
\end{equation}
In the sense of $C(\mathbb{R}_t; H^s_x)$, we have $\Lambda w\big(x,0)=\phi(x)$ on $\mathbb{R}$, and in the sense of $C(\mathbb{R}_x; H^\frac{2s+1}{4}_t)$, we have $\Lambda w(0,t)=f(t)$ on $[0,T]$.  We therefore look to solve $\Lambda w=w$ for some selection of $T$.  By the chain rule and product rule (see below for details), for $\alpha\geq 2$,
\begin{equation} \label{E:142}
\| \Lambda w_1 - \Lambda w_2 \|_{Z} \leq c T^\sigma(\|w_1\|_{Z_T}^{\alpha-1} + \|w_2\|_{Z}^{\alpha-1}) \|w_1-w_2\|_{Z}
\end{equation}
Now choose $T$ small in terms of $\|\phi\|_{H^s(\mathbb{R}^+)}$ and $\|f\|_{H^\frac{2s+1}{4}(\mathbb{R}^+)}$, so that, by \eqref{E:141} and \eqref{E:142}, $\Lambda$ is a contraction, which yields a unique fixed point $u$, which on $[0,T]$ solves the integral equation
  \begin{equation}
\label{E:173}
u(t)=
\begin{aligned}[t]
&e^{it\partial_x^2}\tilde{\phi} + \mathcal{L}(f - e^{i\cdot \partial_x^2}\tilde{\phi}\big|_{x=0}) \\
&- \lambda \mathcal{D}(u|u|^{\alpha-1}) + \lambda \mathcal{L}( \mathcal{D}(u|u|^{\alpha-1})\big|_{x=0})
\end{aligned}
\end{equation}

Let $S$ be the set of all times $T>0$ for which (1) $\exists$ $u\in Z$ such that $u$ solves \eqref{E:173} on $[0,T]$ and (2) for each pair $u_1,u_2\in Z$, such that $u_1$ solves \eqref{E:173} on $[0,T_1]$ with $T_1\leq T$ and $u_2$ solves \eqref{E:173} on $[0,T_2]$ with $T_2\leq T$, we have $u_1=u_2$ on $[0,\min(T_1,T_2)]$.  

We claim that $T$ as given in the above contraction argument is in $S$.  We need only show condition (2).  But the integral equation \eqref{E:173} has a unique solution by the contraction argument in the space $L_{T_m}^qW_x^{s,r}$, where $T_m = \min(T_1,T_2)$, by Lemmas \ref{L:G}\ref{I:Gm}, \ref{L:I}\ref{I:Im}, \ref{L:BF}\ref{I:BFm} and the fact that $\chi_{[0,T_m]}\mathcal{L}g = \chi_{[0,T_m]}\mathcal{L}(\theta_{T_m}g)$, $\chi_{[0,T_m]}\mathcal{D}w = \chi_{[0,T_m]}\mathcal{D}\chi_{[0,T_m]}w$.
Let $T^*=\sup S$.  Define $u^*$ on $[0,T^*)$ by setting, for $t<T^*$, $u^*(t)=u(t)$ for some $u\in Z$ whose existence is given by condition (1); this is well-defined by condition (2). 

Suppose $T^*<\infty$ and $\lim_{t\uparrow T^*} \|u(\cdot, t) \|_{H^s(\mathbb{R}^+)}\neq \infty$.  Then $\exists \; a$ and a sequence $t_n\to T^*$ such that $\|u^*(t_n)\|_{H^s(\mathbb{R}^+)}\leq a$.  By the above existence argument applied at time $t_n$ for $n$ sufficiently large, we obtain a contradiction, as follows.  We shall select $T=t_n$ for $n$ sufficiently large in a moment.  We have, by assumption, $u_1\in Z$ solving the integral equation 
\begin{equation}
\label{E:170}
u_1(t)=
\begin{aligned}[t]
&e^{it\partial_x^2}\tilde{\phi} + \mathcal{L}(f - e^{i\cdot \partial_x^2}\tilde{\phi}\big|_{x=0}) \\
&- \lambda \mathcal{D}(u_1|u_1|^{\alpha-1}) + \lambda \mathcal{L}( \mathcal{D}(u_1|u_1|^{\alpha-1})\big|_{x=0})
\end{aligned}
\end{equation}
on $[0,T]$.  Apply the above existence argument to obtain $u_2\in Z$ solving, on $[T,T+\delta]$, the integral equation
\begin{equation}
\label{E:171}
u_2(t) = 
\begin{aligned}[t]
&e^{i(t-T)\partial_x^2}u(T) + \mathcal{L}^T( f - e^{i(\cdot - T)\partial_x^2}u(T)\big|_{x=0})\\
&- \lambda \mathcal{D}^T(u_2|u_2|^{\alpha-1}) + \lambda \mathcal{L}^T( \mathcal{D}^T(u_2|u_2|^{\alpha-1})\big|_{x=0})
\end{aligned}
\end{equation}
where
$$\mathcal{L}^Tg(t) = \mathcal{}(g(\cdot + T))(t-T) \qquad \mathcal{D}^Tv(t)=\mathcal{D}(v(\cdot+T))(t-T)$$
Since $\delta = \delta(a, \|f\|_{H^\frac{2s+1}{4}(\mathbb{R}^+)})$, we can select $n$ sufficiently large so that $T+\delta=t_n+\delta>T^*$.  Now we show that we can concatenate these two integral equations.  Define $u(t)=u_1(t)$ for $-\infty< t\leq T$ and $u(t)=u_2(t)$ for $T\leq t<+\infty$.  Then clearly $u\in L_t^qW_x^{r,s}\cap C(\mathbb{R}_t;\; H^s_x)$.  Evaluate  \eqref{E:170} at $t=T$, substitute into \eqref{E:171}, and apply the two identities
\begin{equation} \label{E:200}
\begin{aligned}
&\mathcal{L}g(t) = e^{i(t-T)\partial_x^2}\mathcal{L}g(T) - \mathcal{L}^T(g-e^{i(\cdot-T)\partial_x^2}\mathcal{L}g(T)\big|_{x=0})(t) && \text{for }t\geq T \\
&\mathcal{D}v(t) = e^{i(t-T)\partial_x^2}Dv(T) + \mathcal{D}^Tv(t) && \text{for all }t
\end{aligned}
\end{equation}
with $v(t)=-\lambda u|u|^{\alpha-1}(t)$ and $g(t)=f(t)-e^{it\partial_x^2}\tilde{\phi}\big|_{x=0} - \mathcal{D}v(0,t)$ on $0\leq t\leq T+\delta$.  This establishes that $u$ solves, on $[0,T+\delta]$, the integral equation \eqref{E:173}.  
Next, we show that $u \in C(\mathbb{R}_x; \; H^\frac{2s+1}{4})$.  Let $\psi\in C^\infty$ such that $\psi(t)=0$ for $t\leq 0$, $\psi(t)=1$ for $\frac{T}{2}\leq t \leq T+\frac{\delta}{2}$, $\psi(t)=0$ for $t>T+\delta$.  It is clear from the definition of $u$ that $(1-\psi)u \in C(\mathbb{R}_x; \; H_t^\frac{2s+1}{4})$.  Since by \eqref{E:173}
$$ \psi(t)u(t) = 
\begin{aligned}[t]
&\psi(t) e^{it\partial_x^2}\tilde{\phi} + \psi(t) \mathcal{L}( f - \theta_{2(T+\delta)}e^{i\cdot \partial_x^2} \tilde{\phi}\big|_{x=0} )(t) \\
&- \lambda \psi(t) \mathcal{D}( u|u|^{\alpha-1})(t) + \lambda \psi(t) \mathcal{L}( \theta_{2(T+\delta)} \mathcal{D}(u|u|^{\alpha-1})\big|_{x=0})
\end{aligned}
$$  
by Lemmas \ref{L:G}\ref{I:Gt}, \ref{L:I}\ref{I:It}, and \ref{L:BF}\ref{I:BFt} we have $\psi u \in C(\mathbb{R}_x; \; H_t^\frac{2s+1}{4})$.  

Next, we need to verify condition  (2) in the definition of $f$.  Now suppose $u$ is a solution on $[0,T_u]$ with $T_u\leq T+\delta$, and $v$ is a solution on $[0,T_v]$ with $T_v\leq T+\delta$, and suppose $\min(T_u,T_v) \geq T^*$.  Then $u(t)=v(t)$ for all $t\leq T$ (since $T\in S$).  Then, again by \eqref{E:200}, $u$ solves \eqref{E:171} with $u_2$ replaced by $u$, and $v$ solves \eqref{E:171} with $u_2$ replaced by $v$ ($u(T)=v(T)$).  By uniqueness of the fixed point to \eqref{E:171} in $L_{[T,T+\delta]}^qW_x^{s,r}$, we get that $u(t)=v(t)$ on $[T,T+\delta]$.  We have thus established that $\sup S\geq T+\delta>T^*$, which is a contradiction, so in fact $\lim_{t\uparrow T^*} \|u(\cdot, t) \|_{H^s(\mathbb{R}^+)}= \infty$ if $T^*<\infty$.

Now we move on the continuity claim.  Suppose $(\phi,f)$ gives a solution $u$ of \eqref{E:173} on $[0,T^*)$, and consider $(\phi_1,f_1)$ with $\|\phi-\phi_1\|_{H^s(\mathbb{R}^+)} + \|f-f_1\|_{H^\frac{2s+1}{4}(\mathbb{R}^+)}< \delta$.   Fix $T<T^*$.  Let $u_1$ be the solution corresponding to $(\phi_1,f_1)$ on $[0,T_1]$, where $T_1$ is the first time $t$ such that $\|u_1\|_{L_{[0,t]}^qW_x^{s,r}}=2\|u\|_{L_T^qW_x^{s,r}}$.  We claim that $T_1>T$ provided we take $\delta$ sufficiently small.  Indeed, taking the difference of the two integral equations, we find, for $t\leq \min(T_1,T)$ 
$$\| u-u_1 \|_{L_{[0,t]}^q W_x^{s,r}} \leq c\delta + c(\|u\|_{L_T^qW_x^{s,r}} + \|u_1\|_{L_{T_1}^qW_x^{s,r}}) \|u-u_1\|_{L_{[0,t]}^{q_1}W_x^{s,r}}$$
where $q_1<q$, and $c$ depends only upon operator norms.  This gives, by Lemma \ref{L:Gronwall}, 
\begin{equation} \label{E:190}
\|u-u_1\|_{L_{[0,t]}^q W_x^{s,r}} \leq c\delta
\end{equation}
where now $c$ depends on $f$, $\phi$, and $T$.  Now if $T_1<T$, then take $t=T_1$ in \eqref{E:190} and $\delta$ sufficiently small to obtain a contradiction.  The inequality \eqref{E:190} plus estimates on the difference of the integral equations for $u$ and $u_1$ also shows
$$\|u-u_1\|_{C([0,T]; \; H^s_x)} + \|u-u_1\|_{C(\mathbb{R}_x; \; H^\frac{2s+1}{4}(0,T))}  \leq c\delta$$

Now we remark on the proof in the subcritical case (a) for $\frac{1}{2}<s<\frac{3}{2}$.  Let
$$Z= C(\mathbb{R}_t; \; H^s_x) \cap C( \mathbb{R}_x; H^\frac{2s+1}{4}_t)$$
Set $r=2$, $q=\infty$ in the remainder of the argument above.  Do note, however, that to show $f(t)-\theta_{2T}(t)e^{it\partial_x^2}\tilde{\phi}\big|_{x=0} \in  H^\frac{2s+1}{4}(\mathbb{R}_t^+)$, we need to appeal to the compatibility condition $f(0)=\phi(0)$ and Lemma \ref{JK37}.  Also, by Lemma \ref{JK37}, $\theta_{2T}(t) \mathcal{D}(w|w|^{\alpha-1})(0,t) \in H_0^\frac{2s+1}{4}(\mathbb{R}_t^+)$

Now we discuss the critical case (b).  Let $Z=L_t^qW_x^{s,r}$ with $r=\frac{\alpha+1}{1+(\alpha-1)s}$ and $q=\frac{4(\alpha+1)}{(\alpha-1)(1-2s)}$.  The integral equation is
\begin{equation}\label{E:201}
\Lambda w(t) = 
\begin{aligned}[t]
&\theta_T(t)e^{it\partial_x^2}\tilde{\phi} + \theta_T(t)\mathcal{L}(f- \theta_{2T}e^{i\cdot\partial_x^2}\tilde{\phi}\big|_{x=0})(t) \\
&- \lambda \theta_T(t)\mathcal{D}(w|w|^{\alpha-1})(t) + \lambda \theta_T(t) \mathcal{L}(\theta_{2T}\mathcal{D}(w|w|^{\alpha-1})\big|_{x=0} )(t) 
\end{aligned}
\end{equation}

Now, because $q\neq \infty$,  $\|\theta_T(t)e^{it\partial_x^2}\tilde{\phi}\|_{L_t^qW_x^{s,r}}\to 0$ as $T\downarrow 0$ and $\|\theta_T(t)\mathcal{L}(f- \theta_{2T}e^{i\cdot\partial_x^2}\tilde{\phi}\big|_{x=0})(t)\|_{L_t^qW_x^{r,s}}\to 0$ as $T\downarrow 0$.  Therefore, $\exists\; T>0$ such that 
$$\|\theta_T(t)e^{it\partial_x^2}\tilde{\phi}\|_{L_t^qW_x^{s,r}}+ \|\theta_T(t)\mathcal{L}(f- \theta_{2T}e^{i\cdot\partial_x^2}\tilde{\phi}\big|_{x=0})(t)\|_{L_t^qW_x^{s,r}}< \delta$$
which gives
\begin{equation} 
\|\Lambda w \|_{Z} \leq \delta + c \|w\|_{Z}^\alpha
\end{equation}
For $\delta$ sufficiently small, there will be a fixed point in the space $\{ \, w\in Z \, | \, \|w\|_Z < 2\delta\, \}$.  From $\Lambda u =u$, \eqref{E:201} and Lemmas \ref{L:G}\ref{I:Gs}, \ref{L:I}\ref{I:Is}, \ref{L:BF}\ref{I:BFs}, we can recover the bounds in $C(\mathbb{R}_t; \; H_x^s)$, and by Lemmas \ref{L:G}\ref{I:Gt}, \ref{L:I}\ref{I:It}, \ref{L:BF}\ref{I:BFt}, we can recover the bounds in $C(\mathbb{R}_x; \; H_t^\frac{2s+1}{4})$.  Let $T^*$ be the supremum of all existence times with a uniqueness stipulation, as before.  We are not able to show the blowup statement in this case.  Moreover, we also can only establish the continuity assertion for \textit{some} $T<T^*$.

\subsection{Notes on applying the chain and product rule}
We shall apply the chain rule (Lemma \ref{CWchain}) with $w:\mathbb{R}\to \mathbb{C}$ and $F:\mathbb{C}\to \mathbb{C}$ given by $F(w)=|w|^{\alpha-1}w$, for $\alpha\geq 1$.  Then
$$F'(w) = \begin{bmatrix} (\alpha-1)|w|^{\alpha-3}(\text{Re}\,w)^2 + |w|^{\alpha-1} & (\alpha-1) |w|^{\alpha-3}(\text{Re}\,w)(\text{Im}\,w) \\ (\alpha-1) |w|^{\alpha-3}(\text{Re}\,w)(\text{Im}\,w) &  (\alpha-1)|w|^{\alpha-3}(\text{Im}\,w)^2 + |w|^{\alpha-1} \end{bmatrix}$$
and consequently each component of $F'(w)$ is bounded by $|w|^{\alpha-1}$.  Thus
$$ \| D^s |w|^{\alpha-1}w \|_{L_x^{r'}} \leq c\alpha\| |w|^{\alpha-1} \|_{L_x^{r''}} \|D^s w\|_{L_x^r}$$
where $\frac{1}{r''}=\frac{1}{r'}-\frac{1}{r}=1-\frac{2}{r}$ and $\frac{1}{q''}=\frac{1}{q'}-\frac{1}{q}=1-\frac{2}{q}$.  Since $r,q$ have been selected so that$\frac{1}{(\alpha-1)r''}=\frac{1}{r}-s$ and $\frac{1}{(\alpha-1)q''}>\frac{1}{q}$, we have
$$ \| D^s |w|^{\alpha-1}w \|_{L_x^{r'}} \leq c\|D^sw\|_{L_x^r}^\alpha$$

To handle differences, for $w_0,w_1:\mathbb{R}\to \mathbb{C}$, set $w_\theta  = \theta w_1 + (1-\theta )w_0$.  Then 
$$|w_1|^{\alpha-1}w_1-|w_0|^{\alpha-1}w_0= \int_{\theta =0}^1 (\alpha-1)|w_\theta |^{\alpha-3} w_\theta (w_\theta \circ(w_1-w_0)) + |w_\theta |^{\alpha-1}(w_1-w_0)$$
where $z_1\circ z_2 = (\R z_1)(\R z_2) + (\I z_1)(\I z_2)$.  To this, apply $D^s$, and invoke the product rule (Lemma \ref{CWproduct}) and the chain rule (Lemma \ref{CWchain}).
  
\section{Uniqueness: Proof of Prop.\ \ref{P:uniqueness}} \label{S:uniqueness}

We shall begin by establishing uniqueness of a distributional solution with weak traces for the linear problem for $s\geq 0$.  Given two solutions $u_1$, $u_2$, consider the difference $v=u_1-u_2$.  We are thus assuming 
\begin{equation}
\label{E:130}
v\in C([0,T^*); \; L^2(\mathbb{R}^+))\text{ with }v(x,0)=0
\end{equation}
 and 
\begin{equation}
\label{E:131}
\lim_{x\to 0^+} \|v(x,\cdot)\|_{L^2_{(0,T)}} = 0
\end{equation}
Take $T<T^*$.
Let $\theta(t)$ be a nonnegative smooth function supported on $[-2,-1]$ with $\int \theta =1$.  Let $\theta_\delta(t) = \delta^{-1}\theta(\delta^{-1}t)$.  For $\delta,\epsilon>0$, let
\begin{equation}
\label{E:132}v_{\delta,\epsilon}(x,t) = \iint v(y,s)\theta_\delta(x-y)\theta_\epsilon(t-s) \, dy \, ds
\end{equation}
which defines, in the sense of distributions,  $v_{\delta,\epsilon}(x,t)$ a smooth function on  $-\delta<x<+\infty$, $-\epsilon<t<T-2\epsilon$.  Owing to the assumption \eqref{E:130} we can write
$$v_{\delta,\epsilon}(x,t) = \int_s \theta_\epsilon(t-s) \left[ \int_y v(y,s)\theta_\delta(x-y) \, dy \right] \, ds
$$
where the integrals are defined in the usual sense.  From this it follows that
$$\|v_{\delta,\epsilon}(\cdot,t) \|_{L^2(\mathbb{R}^+_x)} \leq \sup_{t+\epsilon\leq s \leq t+2\epsilon} \|v(\cdot,s) \|_{L^2(\mathbb{R}_x^+)}$$
Owing to the assumption \eqref{E:131},  $\exists \; L>0$ such that $\sup_{0<x\leq 2L} \|v(x,\cdot)\|_{L^2_{(0,T)}} \leq 1$.   It follows that, for $x+2\delta<2L$, \eqref{E:132} can be written
$$v_{\delta,\epsilon}(x,t) = \int_y \theta_\delta(x-y) \left[ \int_s v(y,s) \theta_\epsilon(t-s) \, ds \right] \, dy$$
where the integrals are understood in the usual sense, and we also have
\begin{equation} \label{E:134}
\| v_{\delta,\epsilon}(x,\cdot) \|_{L^2(0,T)} \leq \sup_{x+\delta<y<x+2\delta} \| v(y,\cdot) \|_{L^2(\epsilon,T+2\epsilon)}
\end{equation}
Let
$$v_{\epsilon}(x,t) = \int_s \theta_{\epsilon}(t-s) v(x,s) \, ds$$
which is initially understood as defining, for each $t$, a distribution in $x$ on $(0,+\infty)$.  It follows from \eqref{E:130} that it is also, for each $t$, a square integrable function in $x$ with $\|v_\epsilon(\cdot,t)\|_{L^2(\mathbb{R}^+)} \leq \sup_{t+\epsilon<s<t+2\epsilon} \|v(\cdot,s)\|_{L^2(\mathbb{R}^+)}$ and
\begin{equation} \label{E:136}
\lim_{\epsilon \to 0^+} \|v_\epsilon(\cdot,t)-v(\cdot,t)\|_{L^2(\mathbb{R}^+)} =0
\end{equation}
Now we proceed to the calculation.  The identity is
\begin{equation}
\label{E:135}
\int_0^{+\infty} |v_{\delta,\epsilon}(x,T)|^2 \, dx = \int_{x=0}^{+\infty} |v_{\delta,\epsilon}(x,0)|^2 + 2 \text{Im}\, \int_{t=0}^T \partial_x v_{\delta,\epsilon}(0,t) \overline{v_{\delta,\epsilon}(0,t)} \, dt
\end{equation}
Now $\exists \;x_1$ with $0<x_1<L$ such that $\partial_x v_{\delta,\epsilon}(x_1,t) = L^{-1}(v_{\delta,\epsilon}(L,t)-v_{\delta,\epsilon}(0,t))$, by the mean-value theorem.  Again by the mean-value theorem, $\exists \; x_2$ with $0<x_2<x_1$ such that $\partial_x v_{\delta,\epsilon}(x_1,t)-\partial_x v_{\delta,\epsilon}(0,t) = x_1 \partial_x^2 v_{\delta,\epsilon}(x_2,t)$.  Subtracting,
\begin{equation} 
\label{E:133}
\| \partial_x v_{\delta,\epsilon} (0,\cdot) \|_{L^2(0,T)} \leq L \sup_{0\leq y \leq L} \|\partial_x^2v_{\delta,\epsilon}(y,\cdot)\|_{L^2(0,T)} + L^{-1}\sup_{0\leq y \leq L} \| v_{\delta,\epsilon}(y,\cdot)\|_{L^2(0,T)}
\end{equation}
Bounding the terms on the right of this equation, we have
$$ \sup_{0<x<L} \|v_{\delta,\epsilon}(x,\cdot)\|_{L^2(0,T)} \leq \sup_{\delta<x<L+2\delta} \|v(x,\cdot) \|_{L^2_{(\epsilon,T+2\epsilon)}}\leq 1$$
We also have
$$\partial_x^2 v_{\delta,\epsilon}(x,t) = -i\partial_t v_{\delta,\epsilon}(x,t) = i\epsilon^{-1} \iint \theta_\delta(x-y) (\theta')_\epsilon(t-s) v(y,s) \, dy\, ds$$
and thus
$$\sup_{0<x<L} \| \partial_x^2 v_{\delta,\epsilon}(x,\cdot)\|_{L^2(0,T)} \leq \epsilon^{-1} \sup_{\delta<x<L+2\delta} \|v(x,\cdot)\|_{L^2(\epsilon,T+2\epsilon)}\leq \epsilon^{-1}$$
Hence, for fixed $\epsilon>0$ by Cauchy-Schwarz,  bounding by \eqref{E:133} and \eqref{E:134}, we have
$$\int_0^T v_{\delta,\epsilon}(0,t) \overline{\partial_xv_{\delta,\epsilon}(0,t)} \, dt \to 0 \text{ as } \delta \to 0$$
Send $\delta \to 0$ in \eqref{E:135} to get
$$\int_{x=0}^{+\infty} |v_{\epsilon}(x,T)|^2 = \int_{x=0}^{+\infty} |v_{\epsilon}(x,0)|^2 \, dx$$
and then send $\epsilon \to 0$ and use \eqref{E:136}.

Now we prove Prop.\ \ref{P:uniqueness}
\begin{proof}
Suppose $u_1$, $u_2$ are given as in the statement of the proposition, and additionally are smooth and have adequate decay.  Let $v=u_2-u_1$ so that 
$$i\partial_tv+\partial_x^2v + \lambda(|u_2|^{\alpha-1}u_2-|u_1|^{\alpha-1}u_1)$$
and $v(x,0)=0$, $v(0,t)=0$.  Then
\begin{equation} \label{E:191}
\partial_t \int_0^{+\infty} |v|^2 \, dx = 2\text{Re }i\lambda \int_{x=0}^{+\infty} (u_2|u_2|^{\alpha-1}-u_1|u_1|^{\alpha-1})\bar{v} \, dx
\end{equation}
and thus, for any $t>0$,
$$\|v(t)\|_{L_x^2(\mathbb{R}^+)}^2 \leq 2|\lambda| (\|u_1\|_{L_{[0,t]}^\infty L_x^\infty(\mathbb{R}^+)}^{\alpha-1}+\|u_1\|_{L_{[0,t]}^\infty L_x^\infty(\mathbb{R}^+)}^{\alpha-1}) \int_0^t \|v(s)\|_{L_x^2}^2 \, ds$$
By the Sobolev imbedding $H^s(\mathbb{R}^+) \subset L^\infty(\mathbb{R}^+)$ and Gronwall's inequality, $v(t)=0$.  To handle rough $u_1$, $u_2$, mollify $v$ as was done above in the linear case to obtain $v_{\delta,\epsilon}$ so that
$$\partial_t v_{\delta,\epsilon} = i\partial_x^2v_{\delta,\epsilon} +i\lambda(u_2|u_2|^{\alpha-1}-u_1|u_1|^{\alpha-1})_{\delta,\epsilon}$$
Now prove an identity analogous to \eqref{E:135}, estimate as in \eqref{E:191}, and pass to the limit to conclude $v=0$.

\end{proof}

\def\cprime{$'$} \def\polhk#1{\setbox0=\hbox{#1}{\ooalign{\hidewidth
  \lower1.5ex\hbox{`}\hidewidth\crcr\unhbox0}}}
\providecommand{\bysame}{\leavevmode\hbox to3em{\hrulefill}\thinspace}


\begin{thebibliography}{10}

\bibitem{BSZ04}
J.~L. Bona, S.~M. Sun, and B.-Y. Zhang, \emph{Conditional and
  unconditional well-posedness for nonlinear evolution equations}, Advances in
  Differential Equations \textbf{9} (2004), no.~3-4.

\bibitem{BSZ05}
\bysame, \emph{Boundary smoothing and
  well-posedness of the {K}orteweg-de {V}ries equation in a quarter plane}, in
  preparation.

\bibitem{MR2033706}
A.~Boutet~de Monvel, A.~S. Fokas, and D.~Shepelsky, \emph{Analysis of the
  global relation for the nonlinear {S}chr\"odinger equation on the half-line},
  Lett. Math. Phys. \textbf{65} (2003), no.~3, 199--212. 

\bibitem{MR81i:35139}
H.~Br{\'e}zis and T.~Gallouet, \emph{Nonlinear {S}chr\"odinger evolution
  equations}, Nonlinear Anal. \textbf{4} (1980), no.~4, 677--681.
 

\bibitem{MR2001i:35259}
Q.-Y. Bu, \emph{The nonlinear {S}chr\"odinger equation on the semi-infinite
  line}, Chinese Ann. Math. Ser. A \textbf{21} (2000), no.~4, 437--448.
  

\bibitem{MR92d:35268}
R. Carroll and Q.-Y. Bu, \emph{Solution of the forced nonlinear
  {S}chr\"odinger ({NLS}) equation using {PDE} techniques}, Appl. Anal.
  \textbf{41} (1991), no.~1-4, 33--51. 

\bibitem{CW90}
T. Cazenave and F.~B. Weissler, \emph{The {C}auchy problem for the
  critical nonlinear {S}chr\"odinger equation in {$H\sp s$}}, Nonlinear Anal.
  \textbf{14} (1990), no.~10, 807--836. 

\bibitem{CW91}
F.~M. Christ and M.~I. Weinstein, \emph{Dispersion of small amplitude solutions
  of the generalized {K}orteweg-de {V}ries equation}, J. Funct. Anal.
  \textbf{100} (1991), no.~1, 87--109. 

\bibitem{CK02}
J.~E. Colliander and C.~E. Kenig, \emph{The generalized {K}orteweg-de {V}ries
  equation on the half line}, Comm. Partial Differential Equations \textbf{27}
  (2002), no.~11-12, 2187--2266. 

\bibitem{MR2004d:37100}
A.~S. Fokas, \emph{Integrable nonlinear evolution equations on the half-line},
  Comm. Math. Phys. \textbf{230} (2002), no.~1, 1--39. 

\bibitem{F98}
F.~G. Friedlander, ``Introduction to the theory of distributions'', second
  ed., Cambridge University Press, Cambridge, 1998, With additional material by
  M. Joshi. 

\bibitem{H04}
J. Holmer, \emph{The initial-boundary value problem for the {K}orteweg-de
  {V}ries equation}, in preparation.

\bibitem{JK95}
D. Jerison and C.~E. Kenig, \emph{The inhomogeneous {D}irichlet problem
  in {L}ipschitz domains}, J. Funct. Anal. \textbf{130} (1995), no.~1,
  161--219. 

\bibitem{KT98}
M. Keel and T. Tao, \emph{Endpoint {S}trichartz estimates}, Amer. J.
  Math. \textbf{120} (1998), no.~5, 955--980. 

\bibitem{KPV91}
C.~E. Kenig, Gustavo Ponce, and Luis Vega, \emph{Oscillatory integrals and
  regularity of dispersive equations}, Indiana Univ. Math. J. \textbf{40}
  (1991), no.~1, 33--69. 

\bibitem{KPV93}
\bysame, \emph{Small solutions to nonlinear {S}chr\"odinger equations}, Ann.
  Inst. H. Poincar\'e Anal. Non Lin\'eaire \textbf{10} (1993), no.~3, 255--288.
  

\bibitem{St70}
E.~M. Stein, ``Singular integrals and differentiability properties of
  functions'', Princeton Mathematical Series, No. 30, Princeton University
  Press, Princeton, N.J., 1970. 

\bibitem{MR2002d:35196}
W. Strauss and C. Bu, \emph{An inhomogeneous boundary value problem
  for nonlinear {S}chr\"odinger equations}, J. Differential Equations
  \textbf{173} (2001), no.~1, 79--91. 

\bibitem{S77}
R.~S. Strichartz, \emph{Restrictions of {F}ourier transforms to quadratic
  surfaces and decay of solutions of wave equations}, Duke Math. J. \textbf{44}
  (1977), no.~3, 705--714. 

\bibitem{MR92j:35175}
M.~Tsutsumi, \emph{On global solutions to the initial-boundary value problem
  for the nonlinear {S}chr\"odinger equations in exterior domains}, Comm.
  Partial Differential Equations \textbf{16} (1991), no.~6-7, 885--907.
  

\bibitem{MR90h:35224}
M. Tsutsumi, \emph{On smooth solutions to the initial-boundary value
  problem for the nonlinear {S}chr\"odinger equation in two space dimensions},
  Nonlinear Anal. \textbf{13} (1989), no.~9, 1051--1056.

\bibitem{MR85f:35064}
Y. Tsutsumi, \emph{Global solutions of the nonlinear {S}chr\"odinger
  equation in exterior domains}, Comm. Partial Differential Equations
  \textbf{8} (1983), no.~12, 1337--1374. 

\bibitem{MR2001m:35291}
B. Wang, \emph{On the initial-boundary value problems for nonlinear
  {S}chr\"odinger equations}, Adv. Math. (China) \textbf{29} (2000), no.~5,
  421--424. 

\end{thebibliography}
\end{document}